\documentclass{amsart}
\usepackage{amssymb,latexsym}

\theoremstyle{plain}
\newtheorem{theorem}{Theorem}

\newtheorem{remark}{Remark}

\newtheorem{example}[theorem]{Example}
\newtheorem{argument}[theorem]{Argument}

\theoremstyle{definition}
\newtheorem{definition}{Definition}

\theoremstyle{remark}

\numberwithin{equation}{section}

\begin{document}
\title[ Logic With Verbs and its Mathematical Structure]{Logic With Verbs and its Mathematical Structure}
\author{Jun Tanaka}
\address{University of California, Riverside, USA}
\email{juntanaka@math.ucr.edu, junextension@hotmail.com}

\date{Jan, 01, 2010}

\begin{abstract}
The aim of this paper is to introduce the idea of Logic with Verbs and to show its mathematical structure.
\end{abstract}

\maketitle


\section{\textbf{Introduction}}\label{SyloS:1}

In this paper, we introduce the idea of Logic With Verbs as well as its mathematical structure, particularly the negation and the contraposition in Logic With Verbs. Furthermore, we will investigate the soundness of the equivalence between conditional statements (P $\Rightarrow$ Q ) and ($\neg$ P or Q) in Logic With Verbs as well as its Boolean Algebraic Structure. In later sections of this paper, we include observations of the relationship between logic and conversations as well as a discussion for applications of this modern logic and for future research.

\subsection{\textbf{The general form of Logic With Verbs}}
In the paper \cite{verblogic}, the author presented an example of Logic With Verbs as follows:

Premise 1: Tokyo is a part of Japan and Los Angeles is a part of U.S.

Premise 2: Flying is a way of traveling.

$\underline{\text{Premise 3: I flew from Tokyo to Los Angeles.  }}$

Conclusion: I traveled from Japan to U.S.

\subsection{\textbf{The contraposition of Logic With Verbs}}

The negation in Logic With Verbs is as follows:

Premise 1: A hybrid car is a kind of car.

Premise 2: Buying is a way of owning.

$\underline{\text{Premise 3: I have never owned a car}}$

Conclusion: I have never bought a hybrid car.

\subsection{\textbf{The Boolean Algebraic Structure}}

The key idea of this paper is the following argument, which works similarly to that of Classical Logic:

Premise 1: A hybrid car is a kind of car.

$\underline{\text{Premise 2: Buying is a way of owning.}}$

Conclusion: You have never bought a hybrid car or have owned a car.

If Premise 1 and 2 are sound, then the conclusion is sound. This confirms the soundness of the equivalence between conditional statements (P $\Rightarrow$ Q ) and ($\neg$ P or Q) in Logic With Verbs as well as its Boolean Algebraic Structure.

\subsection{\textbf{Abstract and Contents}}
Section 2 is Preliminary for Boolean Algebra and Logic with Verbs. In Section \ref{SyloS:neg}, we will introduce the negation of Logic with Verbs as well as its Boolean Algebraic structure. Starting with partially ordered sets of nouns $\{ N_{i} \}$ and verbs $\{ V_{i} \}$ where the negations $\neg_{n}$ (for nouns) and $\neg_{v}$ (for verbs) satisfy the following properties;
\begin{itemize}
\item $N_{l} < N_{m} \Leftrightarrow \neg_{n} N_{m} < \neg_{n} N_{l}$ \textit{(The law of contrapositive)}
\item $V_{i} < V_{j} \Leftrightarrow \neg_{v} V_{j} < \neg_{v} V_{i}$ \textit{(The law of contrapositive)}
\end{itemize}

We define the composition of a noun and a verb with an operation, which we call Verb Phrase (or simply VP \cite{Chomsky}). Define a binary operator * by * : verb $\times$ noun $\longrightarrow$ Verb Phrase and create partial order $<$ and a negation operator $\neg$ on VPs as follows;

\[
V_{i}* N_{l}  < V_{j}*N_{l} < V_{j}*N_{m} \ \text{and} \
V_{i}* N_{l}  < V_{i}*N_{m} < V_{j}*N_{m}
\]

\[
\neg V_{j}*N_{m}  < \neg V_{j}*N_{l} < \neg V_{i}* N_{l} \ \text{and} \
\neg V_{j}*N_{m}  < \neg V_{i}*N_{m} < \neg V_{i}* N_{l}
\]

In Section \ref{SyloS:1}, we will introduce the mathematical structure of sentences, particularly which have the simplest form as `` Subject Verb Object". Each Verb Prase (Verb*Noun) is expressed on a cartesian product of verb space and noun space, and the negation and the partial order is defined on the cartesian product. Then we will suppose and show several law of Boolean Algebra. We call the collection $\{ V_{a}*N_{b} \}$ of such Verb Phrases (for simplicity VP) with the above properties VP space. For simplicity, in this paper we handle only VP which has the structure as `` verb $+$ one noun". In Section \ref{SyloS:3} we will present the symbolic structure of Logic With Verbs.   In Section \ref{SyloS:6}, \ref{SyloS:7}, and \ref{SyloS:8}, we will discuss some potential methods on how to apply this idea of modern logic to studies in AI communication. In section \ref{SyloS:4}, we will observe the relationship between verbs and nouns. Furthermore, we will discuss how we mutually define verbs and nouns, and also will present a potential application of this modern logic to Fuzzy Set Theory. In section \ref{SyloS:5}, we will introduce a method to relate subjects to verbs and nouns. Recently, a computability of Natural Language is required especially in AI communication
theories. We will introduce several potential approaches on this paper, which
we hope will be a productive contribution to AI in the future.

\subsection{\textbf{The regular form of Logic With Verbs}}\label{SyloS:genform}

In this section, we will relate nouns and verbs from a Set Theoretic view point. Please consider the following three orders of specification;
\begin{itemize}
\item \textit{Orange $<$ fruit $<$ food $($Noun$)$}
\item \textit{Carrot $<$ vegetable $<$ food $($Noun$)$}
\item \textit{Fly $<$ Travel $<$ Move $($Verb$)$}
\end{itemize}

  We will interpret the containments in Set Theory as specificities in order to generalize our usage. A carrot is one kind of vegetable and vegetables are one kind of food. Similarly, to fly is one way to travel, to travel is one way to move. These are merely orders of specificities, and this interpretation of specificity would be more suitable when we apply this Set Theoretic idea to a deductive reasoning as follows.

\begin{itemize}
 \item \textit{I flew from Tokyo to Los Angeles}
\item $\Rightarrow$ \textit{I traveled from Tokyo to Los Angeles} \\  \textit{(By considering flying as a way of traveling)}
\item $\Rightarrow$ \textit{I traveled from Japan to U.S.} \\ \textit{ (By considering Tokyo $<$ Japan and Los Angeles $<$ U.S.)}
\end{itemize}

  Please note that the degree of meaning between the verbs fly, drive, run, and walk will depend on the relative distance to travel from Point A to Point B. Considering the above detailed example, flying is the most suitable way of traveling. This pattern of logic is applicable to the following verbs.

\begin{itemize}
\item \textit{Fly $<$ Travel}
\item \textit{Drive $<$ Travel}
\item \textit{Walk $<$ Travel}
\item \textit{Run $<$ Travel}
\end{itemize}

Just as in Classical Logic, Premise 1 and 2 must be sound. If Premise 3 is factual and the argument is valid, then we say that the conclusion is factual. We use the word sound and factual instead of truth as in Classical Logic because if the statement makes sense for the person or does not, or if the event happened or not is more important in Logic with Verbs than truth. ``True" of Classical Logic works only in the limited world. We rather avoid the long discussion on it in this paper. We will give some examples of modern logic that is presented in the paper \cite{verblogic};

\begin{example}{\textbf{A Regular Form of Logic with Verbs}}

Premise 1: My brother is a lawyer.

Premise 2: Punching is a way of hitting.

$\underline{\text{Premise 3: I punched my brother.   }}$

Conclusion: I hit a lawyer.
\end{example}

\begin{example}{\textbf{A Regular Form of Logic with Verbs}}

Premise 1: A sofa is furniture.

Premise 2: Wiping with a duster is a way of cleaning.

$\underline{\text{Premise 3: I wiped a sofa with a duster.   }}$

Conclusion: I cleaned furniture.
\end{example}

\begin{example}{\textbf{A Regular Form of Logic with Verbs}}\label{example:potato}

Premise 1: A potato is a vegetable.

Premise 2: Baking is a way of cooking.

$\underline{\text{Premise 3: I baked a potato.    }}$

Conclusion: I cooked a vegetable.
\end{example}

 We used a past tense statement for Premise 3 in the above examples since facts are events which happened in the past. Thus, past tense statement is suitable for Premise 3. However, this argument works even with future and present tense sentences as well as sentences with auxiliary verbs as follows:

\begin{example}{\textbf{A Future Tense Form of Logic with Verbs}}

Premise 1: A sofa is furniture.

Premise 2: Wiping with a duster is a way of cleaning.

$\underline{\text{Premise 3: I will wipe a sofa with a duster.   }}$

Conclusion: I will clean furniture.
\end{example}

\subsection{\textbf{The negation in Logic With Verbs}}\label{SyloS:neg}

In this section, we will introduce how to use the negation in Logic With Verbs and express ``if X, and then Y" statements with the negation, and, as well as or; the negation and the expression works similarly to those of Classical Logic.

\begin{example}The negation in Logic With Verbs

\begin{itemize}
 \item \textit{I have never owned a car.}
\item $\Rightarrow$ \textit{I have never bought a car.} \\  \textit{(By considering not owning as a way of not buying.)}
\item $\Rightarrow$ \textit{I have never bought a hybrid car.} \\ \textit{ (By considering $\neg$ car $<$ $\neg$ hybrid car.)}
\end{itemize}
\end{example}

The following two arguments are supposed:

$\text{Buying is a way of Owning} \Leftrightarrow \text{ Not owning is a way of Not buying}$

$\text{A hybrid car is a kind of a car} \Leftrightarrow \text{ Not a car is  Not a hybrid car}$

We have the following examples of negated forms:

\begin{example}{\textbf{A negated form of Logic with Verbs}}\label{def:didnot}

Premise 1: A potato is a vegetable.

Premise 2: Baking is a way of cooking.

$\underline{\text{Premise 3:  I did not cook a vegetable.   }}$

Conclusion: I did not bake a potato.
\end{example}

\begin{example}{\textbf{A past perfect tense negated form of Logic with Verbs}}\label{def:havenever}

Premise 1: A potato is a vegetable.

Premise 2: Baking is a way of cooking.

$\underline{\text{Premise 3:  I have never cooked a vegetable.   }}$

Conclusion: I have never baked a potato.
\end{example}

\begin{example}{\textbf{A past perfect tense negated form; Verb + Two Objects}}\label{ex:TwoObjects}

Premise 1: Tokyo ia a part of Japan.

Premise 2: California is a part of U.S.

Premise 3; Flying is a way of Traveling.

$\underline{\text{Premise 3:  I have never traveled from Japan to U.S.  }}$

Conclusion: I have never flew from Tokyo to L.A.
\end{example}

We consider past perfect tense statements as most suitable tense for Logic with Verbs. Example \ref{def:havenever} may sound more accurate and realistic than Example \ref{def:didnot} because ``have never" shows a experiential fact discussed within the time frame since the person was born before the present, even while the time frame for ``did not" is very vague and it must be implicitly determined by the situation and the communicators. For further discussion on this, please refer to subsection \ref{SyloS:laws}.

\begin{example}{\textbf{logic with Verbs with an intransitive Verb}}

Premise 1: L.A. is a part of California.

Premise 2: Living in X is a way of being to X.

$\underline{\text{Premise 3:  I have never been to California. }}$

Conclusion: I have never lived in L.A.
\end{example}


\section{\textbf{Preliminaries}}\label{SyloS:2}

\subsection{\textbf{Boolean Algebra}}\label{SyloS:boolean}

\begin{definition} {Boolean Algebra}

Boolean algebra provides the operations and the rules for working with the set $\{$0,1$\}$. The complement of an element, denoted with $\neg$, is defined by $\neg$ 0 = 1 and $\neg$ 1 = 0. The Boolean sum, denoted by + or by OR, has the following values:
\[
1 + 1 = 1, \ \ 1 + 0 = 1, \ \ 0 + 1 = 1, \ \ 0 + 0 = 0
\]
The Boolean product, denoted by $\cdot$ or AND, has the following values:
\[
1 \cdot 1 = 1, \ \ 1 \cdot 0 = 0, \ \ 0 \cdot 1 = 0, \ \ 0 \cdot 0 = 0
\]

\end{definition}

\begin{definition} {The Abstract Definition of Boolean Algebra}\label{SyloS:boolean2}

A Boolean Algebra is a set B with two binary operations $\wedge $ and $\vee$, elements 0 and 1, and a unitary operation $\neg$ such that these properties hold for all x, y, and z in B:
\[
\begin{aligned}
 x \ \vee \ 0 = x &    \ \ \ \ \ \ \text{Identity Laws}& \\
 x \ \wedge \ 1 = x &    \ \ \ \  \ \ \text{Identity Laws}&\\
 x \ \vee \ \neg x = 1 &    \ \ \ \ \ \ \text{The law of excluded middle}& \\
 x \ \wedge \ \neg x = 0 &    \ \ \ \  \ \ \text{The law of non-contradiction}&\\
 x \ \vee \ y = y \ \vee \ x  &    \ \ \ \ \ \ \text{Commutative laws}& \\
 x \ \wedge \ y = y \ \wedge \ x &    \ \ \ \  \ \ \text{Commutative laws}&\\
\end{aligned}
\]

\[
\begin{aligned}
 (x \ \vee \ y) \vee \ z = x \ \vee \ (y \vee \ z) &    \ \ \ \ \ \ \text{Associative laws}& \\
 (x \ \wedge \ y) \wedge \ z = x \ \wedge \ (y \wedge \ z) &    \ \ \ \  \ \ \text{Associate laws}&\\
 x \ \vee ( \ y \wedge \ z ) = ( x \ \vee \ y ) \wedge ( x \ \vee \ z) &    \ \ \ \ \ \ \text{Distributive laws}& \\
 x \ \wedge (\ y \vee \ z) = (x \ \wedge \ y) \vee ( x \ \wedge \ z) &    \ \ \ \  \ \ \text{Distributive laws}&\\
\end{aligned}
\]

\end{definition}

\section{\textbf{Mathematical Structure of Logic With Verbs and further discussion}}\label{SyloS:1}
The aim of this paper is to show the mathematical structure of Logic with Verbs. As for Logic with Verbs, please refer to \cite{verblogic}.

\subsection{\textbf{Noun Space and Verb Space}}\label{LWVmath:1}
Let N be a well-defined set with a partial order $\leq_{n}$ and a negation $\neg_{n}$ and be closed under $\neg_{n}$, denoted by (N,$\leq_{n}$,$\neg_{n}$), which satisfies the following property;

For any $N_{l} , N_{m} \in N$
\begin{itemize}
\item $N_{l} \leq_{n} N_{l}$
\item $N_{l} \leq_{n} N_{m}  \ \text{and} \ N_{m} \leq_{n} N_{n} \Rightarrow N_{l} \leq_{n} N_{n}$
\item $N_{l} \leq_{n} N_{m}  \Leftrightarrow \neg_{n} N_{m} \leq_{n} \neg_{n} N_{l}$ \textit{(The law of contrapositive)}
\item $\neg_{n} \neg_{n} N_{m} = N_{m}$ \textit{(The law of double negation)}
\end{itemize}

Similarly to (N,$\leq_{n}$,$\neg_{n}$), we define a well defined space (V,$\leq_{v}$,$\neg_{v}$) with a partial order and a negation and closed under $\neg_{v}$ as follows;
For any $V_{i} , V_{j} \in V$
\begin{itemize}
\item $V_{i} \leq_{v} V_{i}$
\item $V_{l} \leq_{n} V_{m} \ \text{and} \ V_{m} \leq_{n} V_{n} \Rightarrow V_{l} \leq_{n} V_{n}$
\item $V_{i} \leq_{v} V_{j}  \Leftrightarrow \neg_{v} V_{j} \leq_{v} \neg_{v} V_{i}$ \textit{(The law of contrapositive)}
\item  $\neg_{v} \neg_{v} V_{j} = V_{j}$ \textit{(The law of double negation)}
\end{itemize}

(N,$\leq_{n}$,$\neg_{n}$) is called Noun Space and (V,$\leq_{v}$,$\neg_{v}$) is called Verb Space, and for an application of Natural Language Process and Linguistics, N represents a set of nouns and V represents a set of verbs.

\subsection{\textbf{Verb Phrase Space}}\label{LWVmath:2}
Now we will construct a space ($\mathcal{VP}$,$\leq$,$\neg$), called Verb Phrase space, where $\mathcal{VP}$ :=N$\times$V is defined on the cartesian product of Noun Space and Verb Space with the following definition;

(1) $\neg$($V_{i}$,$N_{l}$):=($\neg_{v}$ $V_{i}$,$\neg_{n}$ $N_{l}$)

(2) ($V_{1}$,$N_{1}$) $\leq$ ($V_{2}$,$N_{2}$) if $V_{1}$ $\leq_{v}$ $V_{2}$ and $N_{1}$ $\leq_{n}$ $N_{2}$

$\\$

The law of contrapositive and the law of double negation for Verb Phrase are derived as in the following theorem.
\begin{theorem} Let ($\mathcal{VP}$,$\leq$,$\neg$) be a Verb Phrase space and Let ($V_{i}$,$N_{l}$) $\in$ $\mathcal{VP}$ for all i,l = 1,2,3,$\cdots$.
($V_{1}$,$N_{1}$) $\leq$ ($V_{2}$,$N_{2}$) $\Leftrightarrow$ $\neg$($V_{2}$,$N_{2}$) $\leq$ $\neg$($V_{1}$,$N_{1}$), and $\neg$ $\neg$($V_{i}$,$N_{l}$):= ($V_{i}$,$N_{l}$). Furthermore, every Verb Phrase space is well-defined and closed under $\neg$.
\begin{proof}
The claim follows from the above properties.
\end{proof}
\end{theorem}

\begin{theorem}Let ($\mathcal{VP}$,$\leq$,$\neg$) be a Verb Phrase space and Let ($V_{i}$,$N_{l}$) $\in$ $\mathcal{VP}$ for all i,l = 1,2,3,$\cdots$.
($V_{1}$,$N_{1}$) $\leq$  ($V_{1}$,$N_{2}$) $\leq$ ($V_{2}$,$N_{2}$) and ($V_{1}$,$N_{1}$) $\leq$  ($V_{2}$,$N_{1}$) $\leq$ ($V_{2}$,$N_{2}$)

$\neg$($V_{2}$,$N_{2}$) $\leq$ $\neg$($V_{1}$,$N_{2}$) $\leq$ $\neg$($V_{1}$,$N_{1}$) and $\neg$($V_{2}$,$N_{2}$) $\leq$ $\neg$($V_{2}$,$N_{1}$) $\leq$ $\neg$($V_{1}$,$N_{1}$)

\begin{proof}
Obvious
\end{proof}
\end{theorem}

For simplicity in presentation, ($V_{i}$,$N_{j}$) will be written as $V_{i} \ast N_{j}$ throughout the remainder of this paper. A Verb Phrase attached with a subject at the beginning is called a sentence; For example, I $V_{i} \ast N_{j}$. please note that I $V_{i} \ast N_{j}$ look like a sentence with the subject I. For sentences, $\leq$ may be written with $\Longrightarrow$; please note that $I \ A* E  \Longrightarrow I \ B*F$ is more intuitively clear regarding the flow of the argument.

Furthermore, Min and Max are required to construct Lattice on $\mathcal{VP}$. We can suppose X has done something is Max and X has not done anything is Min. In order to construct Lattice, we suppose the following conditions;

(1) ($V_{i}$,$N_{l}$) $\leq$ ($V_{\text{do}}$,$N_{\text{something}}$) for all i and l.

(2) $\neg$ ($V_{\text{do}}$,$N_{\text{something}}$) $\leq$ $\neg$ ($V_{i}$,$N_{l}$) for all i and l.

X $\neg$ $V_{\text{do}} * N_{\text{something}}$ is supposed to be rendered to X has not done anything.

\subsection{\textbf{Symbolic Structure of Logic With Verbs}}\label{SyloS:3}
In Logic With Verbs, we fix subject and discuss the connection (VP) between verbs and nouns as well as the validity of the flow from one statement to the other. Thus, `` A subject + VP" is called a sentence. Every sentence is factual (strictly) or not factual. For a fixed subject YOU (the readers), let A,B, be verbs where A $\leq$ B and $\leq$ be a partial order. In other words, `` Aing implies Bing" is sound for the readers. Let E,F be nouns where  and E $\leq$ F and $\leq$ be a partial order. In other words, ``E implies F" is sound for the readers. Then we have
\[
A* E  \Longrightarrow B*E \Longrightarrow B*F \ \text{and} \
A* E  \Longrightarrow A*F \Longrightarrow B*F
\]
That means, in terms of sentences,
\[
\text{If ``you A* E" is factual, then ``you B*E" is factual and ``you B*F" is factual} \ \text{ and ``you A*F" is factual.}
\]
\[
\text{If ``you B*E" is factual, then ``you B*F" is factual.}
\]
\[
\text{If ``you A*F" is factual, then ``you B*F" is factual.}
\]

This woks very similarly to Classical Logic but the word ``factual" is used. Logic with Verbs is not made to say ``true" but only to discuss facts.

As for the negation $\neg$,
\[
\neg B*F  \Longrightarrow \neg B*E \Longrightarrow \neg A* E \ \text{and} \
\neg B*F  \Longrightarrow \neg A*F \Longrightarrow \neg A* E
\]
means
\[
\text{If ``you $\neg$ B*F" is factual, then ``you $\neg$ B*E" is factual and ``you $\neg$ A* E" is factual,}
\]
and the rest of arguments are omitted since they would be driven similarly to the previous argument.

By supposing the law of non-contradiction and the law of excluded middle on VPs, we make the following assumptions;

for any subject X, any VP Y*Z
\[
\text{ ``X Y*Z" is factual if and only if ``X $\neg$ Y*Z" is not factual.}
\]
Furthermore,
\[
\text{ ``X $\neg$ Y*Z" is factual if and only if ``X Y*Z" is not factual.}
\]
As for the validity of these two laws, please refer to Section \ref{SyloS:laws}.

\subsection{\textbf{Definition of And, Or in VP Space}}\label{SyloS:andor}
In this section, we will redefine AND as well as OR of written languages (Later we will call linguistic formed sentences). Since linguistic formed sentences like ``I baked potatoes and apples" can not be handled by itself in Logic with Verbs. Thus we need to redefine And as well as Or of linguistic formed sentences in order to transform them into a suitable form of Logic With Verbs. When we say ``I cooked vegetables and fruits", we are not thinking of an intersection of vegetables and fruits as in Classical Logic. We understand and define the sentence ``I cooked vegetables and fruits" as ``I cooked vegetables and I cooked fruits" since it would be more natural to understand the sentence ``I cooked vegetables and fruits" is merely a simplification of combining the two sentences. VP space must be closed under the binary operations AND as well as OR. Thus VP A AND VP B is a VP. Sentence A AND Sentence B is a sentence. So let sentence A and sentence B be VPs with a subject. If ``sentence A and sentence B" is factual, then we say both of sentences are factual. Further, if ``sentence A or sentence B" is factual, then either sentences, possibly both, is factual. We could create a factual table for AND as well as OR between two sentences just as a truth table for AND as well as OR in Classical Logic. Associative Law and Distributive law hold. The proof is driven just as these laws in Classical Logic.

By supposing the law of non-contradiction and the law of excluded middle on VPs, we make the following assumptions;

for any subject X, any VP Y*Z
\[
\text{ ``X Y*Z" is factual if and only if ``X $\neg$ Y*Z" is not factual.}
\]
Furthermore,
\[
\text{ ``X $\neg$ Y*Z" is factual if and only if ``X Y*Z" is not factual.}
\]
As for the validity of these two laws, please refer to Section \ref{SyloS:laws}.

\subsection{\textbf{Linguistics formed sentence V.S. sentences of Logic With Verbs}}\label{SyloS:LinvsLGW}

The above example can be extended to the form as of right and left distributive laws:

For a fixed subject I,
Let A,B,C,D be verbs where A $\leq$ B and C $\leq$ D and $\leq$ be a partial order. Let E,F,G,H be nouns where  and E $\leq$ F and G $\leq$ H and $\leq$ be a partial order. Then and as well as or are redefined between sentences in Logic with Verbs: (sentences written on left-hand side of equal sign are in linguistic form, and sentences written on right hand-side of equal sign are in Logic with Verbs form.)
\begin{itemize}
\item \text{A* (E and G) := A*E AND A*G (Left distributive) For example,} \\
\text{I baked potatoes and apples := ``I baked potatoes" AND ``I baked apples."}
\item \text{(A and C)*E := A*E AND C*E (right distributive) For example,} \\
\text{I baked and ate potatoes := ``I baked potatoes" AND ``I ate potatoes."}
\item \text{A* (E or G) := A*E OR A*G (Left distributive) For example,} \\
\text{I baked potatoes or apples := ``I baked potatoes" OR ``I baked apples."}
\item \text{(A or C)*E := A*E OR C*E (right distributive) For example,} \\
\text{I baked or ate potatoes := ``I baked potatoes" OR ``I ate potatoes."}
\end{itemize}

The following example is rendering from a linguistic formed sentence to a linguistic formed sentence through Logic with Verbs:
A* (E and G) := A*E AND A*G $\Rightarrow$  B*E AND B*G =: B*(E and G)

For example, ``I baked potatoes and apples" = ``I baked potatoes" AND ``I baked apples" $\Rightarrow$ ``I cooked vegetable" AND ``I cooked fruit" $\Rightarrow$ ``I cooked vegetable and fruit."

In other words, if ``I baked potatoes and apples" is factual, then ``I cooked vegetable and fruit" is factual.

\subsection{\textbf{Discussion on the tense of sentences as well as the law of non-contradiction and the law of excluded middle in Natural Language}}\label{SyloS:laws}
In Natural Language, we believe that it is fair to accept the law of excluded middle and the law of non-contradiction for VP for the following reasons. (Please refer to Definition \ref{SyloS:boolean2}) For example, Either sentence A ``I have lived in Tokyo" or sentence B ``I have never lived in Tokyo" must be factual. (The law of excluded middle). In addition, sentence A ``I have lived in Tokyo" and sentence B ``I have never lived in Tokyo" can not be factual at the same time. (The law of non-contradiction). As far as past perfect tense sentences go as ``I have done A" and ``I have never done A", the law of excluded middle and the law of non-contradiction work very well as the readers see in this example. Then how about sentences in the other tense? We will observe present continuous tense sentences; for example, sentence C ``I am driving a car" and sentence D ``I am not driving a car". It would be fair to accept that at a certain moment sentence C or sentence D is factual as well as sentence C and sentence D can not be factual at that same time. As long as the verb is present continuous tense, the sentence describes a motion at a certain moment. At the moment, ``the person is doing A" or ``the person is not doing A" is factual as well as ``the person is doing A" and ``the person is not doing A" can not be factual. As far as we handle past perfect and present continuous formed VP sentences, from the previous discussion, Boolean Algebra in the sense of Definition \ref{SyloS:boolean} is established with And, Or, as well as never. Now we will discuss past tense sentences. We already mentioned that the time frame of past tense sentence is vague so that make it very difficult to handle; for example, sentence E ``I ate an apple" and sentence D ``I did not eat an apple" is both possibly factual with different times. So the time frame need to be a little bit more specified; for example, ``I ate an apple yesterday" and ``I did not eat an apple yesterday" can not be factual as well as ``I ate an apple yesterday" or ``I did not eat an apple yesterday" must be factual. We understand that some fuzziness remain in past tense case. The key observation here is that we can establish Boolean Algebra on sentences by specifying the time frame. For sentences in the future tense, sentence ``I will do X" is more less a plan or a thought. We could handle future tense sentences in Logic with Verbs just as the other tense forms, however we rather say the sentence is ``a plan" instead of ``factual".

For further discussion, we would like to mention the followings; this modern logic is not made to handle sentences which express emotion such as ``I am missing her but at the same time I am not missing her." The author personally understands such a moment, however we know that Logic with Verbs does not work properly for most of such literal sentences. Further investigation is required to improve our logic so that we can handle such literal sentences.

\subsection{\textbf{A different expression of conditional sentences}}\label{SyloS:applinegation}
In this section, we show that a conditional ``if and then" sentence can be expressed with OR as well as Negation, supporting the law of non-contradiction and the law of excluded middle. This expression make the structure of Logic With Verbs similar to that of Classical Logic.

\begin{example}{\textbf{Another expression of a conditional sentence}}

Premise 1: A potato is a vegetable.

$\underline{\text{Premise 2: Baking is a way of cooking.  }}$

Conclusion: ``I have never baked a potato" OR ``I have cooked a vegetable."

If Premise 1 and 2 are sound, either sentence A ``I have baked a potato" or sentence B ``I have never baked a potato" must be factual (supposing the law of excluded middle and the law of non-contradiction). If sentence A ``I have baked a potato" is factual, then ``I have cooked a vegetable" must be factual by Premise 1 and 2. Thus we obtain above Conclusion.
\begin{remark}
The conditional sentence ``If I have baked a potato, then I have cooked a vegetable." must be deduced from the above conclusion just as in Classical Logic. From this observation, we obtain the conclusion of the structure of Logic with Verbs; `` $A* E  \Longrightarrow B*F$ " is equivalent to `` $\neg A* E$ OR $B*F$ " . (The symbols are inherited from the section \ref{SyloS:andor} )
\end{remark}
\end{example}

\subsection{\textbf{Conclusion}}\label{SyloS:sec2conclusion}
Identity Law of definition \ref{SyloS:boolean2} is satisfied if 1 is supposed a sentence which is factual and 0 is supposed a sentence which is not factual. From all of the above argument, every law in Boolean Algebra is established.


\section{\textbf{Second Order Logic With Verbs}}
In the previous section, we presented Boolean Algebraic structure of Logic With Verbs which show logical argument flow by sentences consisting of verb and noun. In this section, we will investigate logical argument by sentences of past perfect tense and past tense by using quantifiers. Particularly, the purpose of this section is to analyze second order Logic expression for the following arguments;

\begin{example}{\textbf{The regular form of Logic With Verbs}}\label{Rform:H1}

Premise 1: A laptop computer is a kind of computer.

Premise 2: Buying X (for oneself) is a way of owning X.

$\underline{\text{Premise 3: I have bought a laptop computer}}$

Conclusion: I have owned a computer.

\end{example}

For simplicity, Buying X (for oneself) will be written as buying X throughout the remainder of the section.

\begin{example}{\textbf{The contraposition of Logic With Verbs}}\label{conexa:H1}

The negation in Logic With Verbs is as follows:

Premise 1: A laptop computer is a kind of computer.

Premise 2: Buying X is a way of owning X.

$\underline{\text{Premise 3: I have never owned a computer}}$

Conclusion: I have never bought a laptop computer.

\end{example}

\begin{example}{\textbf{The Boolean Algebraic Structure}}

The key idea of this section is the following argument, which works similarly to that of Classical Logic:

Premise 1: A laptop computer is a kind of computer.

$\underline{\text{Premise 2: Buying X is a way of owning X.}}$

Conclusion: You have never bought a laptop computer or have owned a computer.

\end{example}

\begin{example}{\textbf{Temporality}}

Premise 1: A laptop computer is a kind of computer.

Premise 2: Buying X is a way of owning X.

$\underline{\text{Premise 3: I bought a laptop computer two years ago}}$

Conclusion: I have owned a computer.

\end{example}

\subsection{\textbf{Definition and discussion}}
We interpret statements of past perfect tense as that there is an experience or a time of Ving N, or we could interpret the action V exists at a certain time t. In either way, the interpretation would lead to the same logical expression.  For example, in this section ``I have eaten curry" is interpreted as ``there was a time of eating curry in my life".

Now we suppose the law of non-contradiction and the law of excluded middle, restricted to the time frame as in  \cite{verblogic}.

Thus, by inheriting the notion of Logic With Verbs \cite{verblogic} and with the above interpretation, we will define the statement ``I have $V_{1}*N_{1}$" in a logical manner as the following; For a fixed subject I,
$\exists$ a time t $\in$ [$t_{1}$,$t_{2}$] such that ``I $V_{1,t}*N_{1}$" is factual where [$t_{1}$,$t_{2}$] is the time period when the person of the statement is living. If the person is living, $t_{2}$ is now.

By following the main idea of Logic with Verbs, ``I $V_{1,t}*N_{1}$" is factual implies that ``I $V_{2,t}*N_{2}$" is factual where $V_{1,t} \Rightarrow V_{2,t}$ for all t and $N_{1} \Rightarrow N_{2}$.

From all of the above argument, we render from a linguistic sentence to a logical sentence;

I have bought a laptop computer $\xrightarrow{\mathrm{render}} $ $\exists$ a time t $\in$ [$t_{1}$,$t_{2}$] such that ``I $V_{\text{buy},t}*N_{\text{laptop computer}}$" is factual.

I have not bought a laptop computer $\xrightarrow{\mathrm{render}} $ $\neg \exists$ a time t $\in$ [$t_{1}$,$t_{2}$] such that ``I $V_{\text{buy},t}*N_{\text{laptop computer}}$" is factual. This negation operates as $\forall$ time t $\in$ [$t_{1}$,$t_{2}$] such that ``I $\neg V_{\text{buy},t}*N_{\text{laptop computer}}$" is factual.

Certainly, we can inverse the render from a logical sentence to a linguistic sentence. The inverse is called inverse-render and denoted by $\xrightarrow{\mathrm{invese}}$.

\begin{argument}{\textbf{Logical Argument of Example \ref{Rform:H1}}}

Let's suppose the three following premises;

Premise 1: A laptop computer is a kind of computer. let $N_{\text{laptop computer}}$ be a laptop computer and $ N_{\text{computer}}$ be a computer.

Premise 2: Buying is a way of owning. Let $V_{\text{buy}}$ be buying and $V_{\text{own}}$ be owning.

$\underline{\text{Premise 3: I have bought a laptop computer}}$

$\xrightarrow{\mathrm{render}} $ $\exists$ a time t $\in$ [$t_{1}$,$t_{2}$] such that ``I $V_{\text{buy},t}*N_{\text{laptop computer}}$" is factual.
Hence, $\exists$ a time t $\in$ [$t_{1}$,$t_{2}$] such that ``I $V_{\text{own},t}*N_{\text{computer}}$" is factual.

$\xrightarrow{\mathrm{invese}}$

Conclusion: I have owned a computer.
\end{argument}

\subsection{\textbf{Definition and discussion for negation}}

The negation of the previous statement is ``I have not $V_{2}*N_{2}$" and that is written in a logical manner as follows;

Similarly to the regular second order logic, it is supposed that $\neg$ $\exists$ a time t $\in$ [$t_{1}$,$t_{2}$] such that ``I $V_{2,t}*N_{2}$" is factual $\Leftrightarrow$ $\forall$ time t $\in$ [$t_{1}$,$t_{2}$],  ``I $\neg V_{2,t}*N_{2}$" is factual.

By inheriting the notion of Logic With Verbs \cite{verblogic} ``I $\neg V_{2,t}*N_{2}$" is factual implies ``I $\neg V_{1,t}*N_{1}$" is factual.

Thus, $\forall$ time t $\in$ [$t_{1}$,$t_{2}$],  ``I $\neg V_{1,t}*N_{1}$" is factual.

From all of the above argument, the argument of example \ref{conexa:H1} is constructed as the following example;

\begin{argument}\label{negation:H1}

We suppose the time just as the previous example.

Premise 1: A laptop computer is a kind of computer.

Premise 2: Buying is a way of owning.

$\underline{\text{Premise 3: I have never owned a computer}}$

$\xrightarrow{\mathrm{render}} $ $\neg$ $\exists$ a time t $\in$ [$t_{1}$,$t_{2}$] such that ``I $V_{\text{own},t}*N_{\text{computer}}$" is factual.

$\Rightarrow$ $\forall$ time t $\in$ [$t_{1}$,$t_{2}$],  ``I $\neg V_{\text{own},t}*N_{\text{computer}}$" is factual.

$\Rightarrow$ $\forall$ time t $\in$ [$t_{1}$,$t_{2}$],  ``I $\neg V_{\text{buy},t}*N_{\text{laptop computer}}$" is factual.

$\Rightarrow$ $\neg \exists$ a time t $\in$ [$t_{1}$,$t_{2}$] such that ``I $V_{\text{buy},t}*N_{\text{laptop computer}}$" is factual.

$\xrightarrow{\mathrm{invese}}$

Conclusion: I have never bought a laptop computer.

\end{argument}

\subsection{\textbf{Expression and discussion for Second Order Logic}}
$\\$

\begin{argument}\label{Or:H1}$\\$

Premise 1: A laptop computer is a kind of computer. let $N_{\text{laptop computer}}$ be a laptop computer and $N_{\text{computer}}$ be a computer.

Premise 2: Buying is a way of owning. Let $V_{\text{buy}}$ be buying and $V_{\text{own}}$ be owning.

$\underline{\text{Assumption 1: Either ``I have owned a computer" or ``I have not owned a computer" is factual}}$
$\xrightarrow{\mathrm{render}} $ ``$\exists$ a time t $\in$ [$t_{1}$,$t_{2}$] such that ``X $V_{\text{own},t}*N_{\text{computer}}$" is factual." or ``$\neg \exists$ a time t $\in$ [$t_{1}$,$t_{2}$] such that ``X $V_{\text{own},t}*N_{\text{computer}}$" is factual.

$\Rightarrow$ by Example \ref{negation:H1}, ``$\exists$ a time t $\in$ [$t_{1}$,$t_{2}$] such that ``X $V_{\text{own},t}*N_{\text{computer}}$" is factual." or ``$\neg \exists$ a time t $\in$ [$t_{1}$,$t_{2}$] such that ``X $V_{\text{buy},t}*N_{\text{laptop computer}}$" is factual.

$\xrightarrow{\mathrm{inverse}}$

Conclusion: You have owned a computer or you have not bought a laptop computer.
\end{argument}

We will extend the above argument to each person's world.
\begin{argument}\label{personalOr:H1}
We suppose the list $\{X_i\}$ of people where each $X_i$ represent each person, and each person accept either Premise N or not. Assumption: If $X_i$ accept Premise N and M, then $X_i$ must accept the conclusion generated by the presented logical argument. Each sentence $X_i$ $V_m * N_l$ must be either factual or not factual just as in the previous sections. we call the collection $\mathfrak{W}$ of sentences the world.

Premise 1: A laptop computer is a kind of computer. let $N_{\text{laptop computer}}$ be a laptop computer and $N_{\text{computer}}$ be a computer.

$\underline{\text{Premise 2: Buying is a way of owning. Let $V_{\text{buy}}$ be buying and $V_{\text{own}}$ be owning.}}$

$\xrightarrow{\mathrm{render}} $ $\forall$ person $X_i$ who accept premise 1 and premise 2, ``$\exists$ a time t $\in$ [$t_{1}$,$t_{2}$] such that ``$X_i$ $V_{\text{own},t}*N_{\text{computer}}$" is factual." or ``$\neg \exists$ a time t $\in$ [$t_{1}$,$t_{2}$] such that ``$X_i$ $V_{\text{buy},t}*N_{\text{laptop computer}}$" is factual.

$\xrightarrow{\mathrm{inverse}}$
Conclusion: If person $X_i$ accept Premise 1 and 2, then $X_i$ have owned a computer or you have not bought a laptop computer.

\end{argument}
By the law of double negation, the following statement is going to be equivalent to the above statement; $\neg \exists$ person  $X_i$ who accept premise 1 and premise 2, ``$\neg \exists$ a time t $\in$ [$t_{1}$,$t_{2}$] such that ``$X_i$ $V_{\text{own},t}*N_{\text{computer}}$" is factual." and ``$\exists$ a time t $\in$ [$t_{1}$,$t_{2}$] such that `` $X_i$ $V_{\text{buy},t}*N_{\text{laptop computer}}$" is factual.

Herein we presented Second Order Logic with Verbs.

\begin{argument}{\textbf{Temporality}}\label{temporality:H1}
We will use the same notation for nouns and verbs as in the previous section.

Premise 1: A laptop computer is a kind of computer.

Premise 2: Buying is a way of owning.

$\underline{\text{Premise 3: I bought a laptop computer two years ago}}$

$\xrightarrow{\mathrm{render}}$

$\exists$ t $\in$ [$t_{3}$,$t_{4}$] such that ``I $V_{\text{buy},t}*N_{\text{laptop computer}}$" is factual where [$t_{3}$,$t_{4}$] represents the time of two years ago.

By considering [$t_{3}$,$t_{4}$] $\subset$ [$t_{1}$,$t_{2}$], $\exists$ a time t $\in$ [$t_{1}$,$t_{2}$] such that ``I $V_{\text{buy},t}*N_{\text{laptop computer}}$" is factual.

$\exists$ a time t $\in$ [$t_{1}$,$t_{2}$] such that ``I $V_{\text{own},t}*N_{\text{computer}}$" is factual.

$\xrightarrow{\mathrm{invese}}$

I have owned a computer.

Similarly, one can prove that I have not owned a computer $\Rightarrow$ I did not buy a laptop computer two years ago.
$\neg \exists$ a time t $\in$ [$t_{1}$,$t_{2}$] such that ``I $V_{\text{own},t}*N_{\text{computer}}$" is factual.
$\neg \exists$ a time t $\in$ [$t_{1}$,$t_{2}$] such that ``I $V_{\text{buy},t}*N_{\text{laptop computer}}$" is factual.
Thus $\neg \exists$ t $\in$ [$t_{3}$,$t_{4}$] such that ``I $V_{\text{buy},t}*N_{\text{laptop computer}}$" is factual.
\end{argument}

Furthermore, from the above argument \ref{temporality:H1} we could extend Argument \ref{personalOr:H1} to a general case as follows.

$\forall$ person $X_i$ who accept premise 1 and premise 2, ``$\exists$ a time t $\in$ [$t_{1}$,$t_{2}$] such that ``$X_i$ $V_{\text{own},t}*N_{\text{computer}}$" is factual." or ``$\neg \exists$ a time t $\in$ [$t_{3}$,$t_{4}$] such that ``$X_i$ $V_{\text{buy},t}*N_{\text{laptop computer}}$" is factual. "
$\xrightarrow{\mathrm{inverse}}$
If a person $X_i$ accept Premise 1 and Premise 2, $X_i$ did not buy a laptop computer two years ago or have owned a computer.

\section{\textbf{Observation on how to Apply This Modern Logic to AI}}\label{SyloS:6}

\subsection{\textbf{Questions for more detail information in conversation}}\label{SyloS:6question}
In this section, we will compare daily conversations with this presented modern logic. Our conversation never flows as Examples shown in Section 2. However, we believe that the structure of the presented modern logic is necessary and applicable to AI communication. We do not need to give the most detailed information in our conversations, thus we provide only sufficient information or only a part that he or she would like to emphasize. Then the listener may ask the speaker for more information if he is interested in more detail. I will give one example of a conversation which distinguishes flow of the presented modern logic.

Person A: ``I traveled to U.S."

Person B: ``Where in U.S. did you travel?"

Person A: ``California"

Person B: ``Where did you fly from?"

Person A: ``I flew from Tokyo"

The above conversations sound more natural than the examples presented in Section 2. Regular conversations typically go from a general statement to a more specific statement, depending on how much information is needed or how much interest is showed in, even while the logical statement flows from the specific statement to a more general statement. In order to make AI communicate ``humanistically'', we suggest generating the most specific statement for each fact beforehand, and then we must make it general enough to ``humanize conversations''. In other words, we need some filtering on generated statements before the output of a statement.

\subsection{\textbf{An Application from the Observation of the Previous Section}}\label{SyloS:7}

Here is a systematized application for more natural conversations from the observation of Natural Language as  shown the previous section:

Premise 1: a house is a kind of a property

Premise 2: California is a part of U.S.

Premise 3: buying X (for myself) is a way of owning X.

Premise 4: I will buy a house in California. (a fact related to the above premises)

We will generate the below seven conclusions out of the four premises.

Conclusion 1: I will buy a house in U.S.

Conclusion 2: I will buy a property in California.

Conclusion 3: I will buy a property in U.S.

Conclusion 4: I will own a house in California.

Conclusion 5: I will own a house in U.S.

Conclusion 6: I will own a property in California.

Conclusion 7: I will own a property in U.S.

In order to make this logic conversational, we need to reverse the pattern that is usually seen in logic. We will demonstrate to generate a conversation between a computer program and a person and let Person A be a computer.

We call HOW, WHICH PART, WHAT KIND question operators which reverse A $\leq$ B. For example, if A $\leq$ B which means A is a kind of B, WHICH KIND * B = A, ``WHICH KIND of property will you buy in California?", the answer is ``I will buy a house in California." (WHICH KIND * property $\Rightarrow$ house.)

Person A ``I will own property in U.S."

Person B ``Which part of U.S. will you own property?"

Person A ``I will own a property in California" (WHICH PART * (own * property * U.S.) $\Rightarrow$ own * property * California)

Person B ``How will you own property in California?"

Person A ``I will buy property in California"n (HOW*(own*property*California)$\Rightarrow$ buy*property*California)

Person B ``Which kind of property will you own in California?"

Person A ``I will buy a house in California" (WHICH KIND*(buy*property*California)$\Rightarrow$ buy*house*California)

If Premises 1 to 4 are input beforehand in a program, it systematically generate correspondences just as above.

\subsection{\textbf{``If And Then" sentence In Logic With Verbs}}\label{SyloS:8}

We will introduce an extension of the application from the previous section, which shows how to handle ``if and then" sentences in Logic with Verb. In addition to premises 1 to 4 in the previous section, we will add one more premise as follows:

Premise 5: if I get this job, I will buy a house in California.

then it implies all of the seven following conclusions.

Conclusion 1': If I get this job, I will buy a house in U.S.

Conclusion 2': If I get this job, I will buy a property in California.

Conclusion 3': If I get this job, I will buy a property in U.S.

Conclusion 4': If I get this job, I will own a house in California.

Conclusion 5': If I get this job, I will own a house in U.S.

Conclusion 6': If I get this job, I will own a property in California.

Conclusion 7': If I get this job, I will own a property in U.S.

$\\$

In the following two subsections, we include an direction for future research, mainly regarding the relation between verbs and nouns as well as a way to handle subjects.
\subsection{\textbf{Recursive Definition of Nouns and Verbs}}\label{SyloS:4}

There are some pairs of verbs and nouns which are defined recursively as a pair; We call such recursive definition N-V isomorphism. In this section, we will show how nouns and verbs should be related through a fuzzy set theoretic view. Some examples of N-V isomorphism as follow;

(1) Food is something you eat. Something you eat is most likely food. 

(2) A Beverage is something you drink. Something you drink is most likely a beverage.

(3) Something you ride on is a vehicle. A vehicle is something you ride on.

(4) Something you draw is a drawing. A drawing is something you draw.

(5) Something you sing is a song. A song is something you sing.

Eat and food are N-V isomorphic, and bread is food. Thus I can eat bread, and the statement ``I can eat bread" is sound, (showing possibility). Now, we will show that N-V isomorphism is used to show the degree of possibility with fuzzy sets; Seaweed is food but if ``I" is American, Seaweed is not very familiar as food. Thus the characteristic value of Seaweed as food must be low. Let's say 0.1. Then the statement ``I can eat Seaweed" should be sound, but the statement ``I rarely eat Seaweed" or ``I am less likely to eat Seaweed" are more appropriate. Now some connections between N-V isomorphism and fuzzy sets are apparent.

So let's suppose the characteristic value of chicken as food is 0.95. ``I often eat chicken" must be appropriate. We could let the range of characteristic values between 1-0.7 be ``often", 0.7-0.4 be ``more or less", 0.4-0.2 be ``less likely", 0.2-0.05 be ``rarely", 0.05-0 ``never". Next we can create a Fuzzy Set Theoretic statement such as ``I often eat pizza", ``I rarely eat deer meat", and ``I never eat a book" by following the method of Zadeh. \cite{Zad,Zad4}. 

\subsection{\textbf{Conditional Logic; How to deal with subjects}}\label{SyloS:5}
In this section, we will present one possibility on how to handle subjects. By using the idea presented in the previous subsection, each person has a different value for classification of each object. In this interpretation, subjects affect and control the degree of possibility for doing X. In the previous section, we mentioned ``I rarely eat Sea Weed" or ``I am less likely to eat Seaweed" if ``I" is American. If ``I" is Japanese, ``I sometimes eat Sea Weed" or ``I often eat Seaweed" must be appropriate. Thus, depending on the subject, the degree of possibility of the combination (Verb*noun) must vary.

\section{\textbf{Conclusion and Observation}}\label{SyloS:conclusion}
In the entire Section \ref{SyloS:2}, we tackled systematic expression of Linguistics and showed Main Boolean Algebraic Structure of sentences in Natural Language. We believe that this Modern Logic Theory is articulated with Modern Set \cite{modern} , and it would help us to bring systematic expression of languages closer to the level of sophistication of human conversations. I also strongly believe that this new logic system could open up a new branch of Artificial Intelligence. This Verb Phrase Logic theory is made only for a specific person and tense. However we will need a integrated logic of all. Further investigation in logic and linguistics are required to improve the systematic expression of our rational thought, which in turn is necessary in creating a communicative Artificial Intelligence. I dream of the day when we can create real AI.

\section{\textbf{Acknowledgement}}
 The author would like to thank all professors who gave him very professional advice and suggestions, which he truly believes improved the preparation of this paper. The author is inspired by Philosophers M.Heidegger, L.Wittgenstein, and N. Chomsky. The author gratefully acknowledges Jay Tomioka for his inspiration and T. Hughes for his editorial assistance.



\end{document}